\documentclass[11pt,1pt,nopreprintline]{elsarticle}

\usepackage{hyperref}
\usepackage{amsmath,amssymb,amsopn,amsfonts,mathrsfs,amsbsy,amscd}
\usepackage{longtable}
\usepackage{color}
\usepackage{graphicx}
\usepackage[all,cmtip]{xy}
\usepackage{amsthm}
\usepackage{mathtools}
\usepackage{extarrows}
\theoremstyle{plain}
\newtheorem{theorem}{Theorem}[section]
\newtheorem{proposition}[theorem]{Proposition}
\newtheorem{lemma}[theorem]{Lemma}

\theoremstyle{definition}
\newtheorem{definition}{Definition}[section]
\newtheorem{example}[theorem]{Example}

\newtheorem{remark}[theorem]{Remark}

\newcommand{\id}{\operatorname{Id}}

\numberwithin{equation}{section}

\newcommand{\End}{\operatorname{End}}
\newcommand{\Aut}{\operatorname{Aut}}

\date{}

\begin{document}

\begin{frontmatter}

\title{\textbf{A Born--Poisson structure on the tangent bundle
		of a Koszul--Vinberg manifold}}

\author[l1]{Abdelhak Abouqateb}
\ead{a.abouqateb@uca.ac.ma}
\author[l1,l2]{Charif Bourzik}
\ead{bourzikcharif@gmail.com}

\address[l1]{MAIST Laboratory, Faculty of Sciences and Technologies, Cadi Ayyad University, B.P. 549 Gueliz, 40000 Marrakesh, Morocco.}
\address[l2]{Higher Institute of Nursing and Health Techniques, Marrakesh, Morocco.}

\begin{abstract}
A Koszul--Vinberg (KV) manifold is a manifold $M$ equipped with a
flat torsion-free connection $\nabla$ and a symmetric bivector field
$h$ satisfying a generalized Codazzi equation. When $h$ is non-degenerate, we recover the well-known notion of a pseudo-Hessian manifold. In this paper, we define the notion of Born-Poisson manifolds, consisting of a Poisson bivector and a compatible
complex product structure. We show that the tangent bundle $TM$ of a KV manifold $M$ carries a canonical Born--Poisson (BP) structure. We prove that two KV manifolds are isomorphic if and
only if their canonical tangent BP manifolds are isomorphic.
Moreover, every such BP isomorphism is uniquely a tangent lift of
a KV isomorphism followed by translation by a parallel vector field.
Finally, we characterize the canonical KV Poisson lifts by Euler
homogeneity and compatibility with the canonical complex structure. 
\\
\noindent{\bf Mathematics Subject Classification 2020:} Primary: 53D05, 53C15, 53A15, 53D17; Secondary:  53C55, 53C50.\\
{\bf Keywords:}  Born--Poisson structure; Koszul--Vinberg manifold; Poisson tensor, Tangent bundle; Euler homogeneity.  
\end{abstract}

\end{frontmatter}	
	
\section{Introduction}\label{sec1}
An \emph{affine manifold} is a manifold \(M\) endowed with a flat torsion-free connection \(\nabla\). In a local affine coordinate system \((x_1,\ldots,x_n)\), one has
$\nabla \partial_{x_i}=0$.

Let \(h\) be a symmetric bivector field on an affine manifold \((M,\nabla)\), with associated bundle map \(h_\sharp:T^*M\to TM\) defined by
$
\beta(h_\sharp\alpha)=h(\alpha,\beta)$. 
The tensor \(h\) is called a \emph{Koszul--Vinberg} (KV) bivector if it satisfies the \emph{generalized Codazzi equation}
\begin{equation}\label{codazzi}
	(\nabla_{h_\sharp\alpha}h)(\beta,\gamma)
	=(\nabla_{h_\sharp\beta}h)(\alpha,\gamma),
	\qquad \forall\alpha,\beta,\gamma\in\Omega^1(M).
\end{equation}
The triple \((M,\nabla,h)\) is then called a \emph{KV manifold}. In affine coordinates, equation~\eqref{codazzi} is equivalent to
\begin{equation}\label{eq2}
	\sum_{l=1}^n
	\left(h_{il}\frac{\partial h_{jk}}{\partial x_l}
	-h_{jl}\frac{\partial h_{ik}}{\partial x_l}\right)=0,
	\qquad \forall i,j,k=1,\ldots,n,
\end{equation}
where \(h_{ij}=h(dx_i,dx_j)\).

KV geometry extends pseudo-Hessian geometry as follows. If \(h\) is nondegenerate, then its inverse is a pseudo-Hessian metric, and \((M,\nabla,h^{-1})\) is a pseudo-Hessian manifold \cite{Shima}. It should be noted (cf. \cite{ABB2020}) that the characteristic distribution \(h_{\sharp}(T^{*}M)\subset TM\) of a KV manifold \((M,\nabla,h)\) is integrable, and its maximal integral submanifolds are the \emph{pseudo-Hessian leaves} of \((M,\nabla,h)\). They foliate \(M\), and each leaf \(L\hookrightarrow M\) is an affine submanifold and carries a canonical pseudo-Hessian metric \(g_L\) characterized by
\[
g_L\bigl(h_{\sharp}\alpha,h_{\sharp}\beta\bigr)=h(\alpha,\beta),
\qquad \alpha,\beta\in T^{*}M.
\]
For further developments in KV geometry, see \cite{BB2019,ABB2022,ABB2023}.

To any KV structure $(\nabla,h)$ on $M$ we can associate a Poisson tensor field $\Pi$ on the tangent bundle \(TM\) \cite{ABB2020}. To explain this, first recall that if \((x_1,\dots,x_n)\) is an affine coordinate system on an open set \(U\subset M\), then the induced tangent coordinates \((x_1,\dots,x_n,y_1,\dots,y_n)\) on \(TU\) are given by
\[
x_i(u):=x_i(p(u)),
\qquad
y_i(u):=(dx_i)_{p(u)}(u),
\]
for \(u\in TU\), where \(p\colon TM\to M\) is the canonical projection. In these tangent coordinates, the expression of the tangent Poisson lift \(\Pi\) on \(TM\) is given by
\begin{equation}\label{eq:intro-lift}
	\Pi=\sum_{i,j=1}^n h_{ij}
	\partial_{x_i}\wedge\partial_{y_j}.
\end{equation}

A \emph{KV isomorphism} is a diffeomorphism that preserves both the connection and the KV bivector. The induced map on tangent bundles preserves the associated Poisson bivector; hence, isomorphic KV manifolds have isomorphic associated tangent Poisson manifolds. However, non-isomorphic KV manifolds may still have isomorphic associated tangent Poisson manifolds. Indeed, Examples~\ref{ex:same-connection} and~\ref{ex:different-connections} show that the Poisson lift alone may fail to distinguish either the symmetric bivector or the affine connection.

\begin{example}[Same connection, different tensors]
	\label{ex:same-connection}
	Let $\nabla$ be the standard connection on $\mathbb R$, and set
	\[
	h_0=\partial_x\otimes\partial_x,\qquad
	h_1=e^x\partial_x\otimes\partial_x.
	\]
	An affine diffeomorphism of $\mathbb R$ has the form $f(x)=ax+b$,
	$a\ne0$. The equality $f_*h_0=h_1$ would require
	$a^2=e^{ax+b}$ for every $x$, which is impossible. Hence the KV
	manifolds $(\mathbb R,\nabla,h_0)$ and $(\mathbb R,\nabla,h_1)$ are not isomorphic. Nevertheless, the diffeomorphism $F:\mathbb R^2\to \mathbb R^2$ given by $F(x,y)=(x,e^x y)$	intertwines their Poisson lifts. Indeed, $F_*(\Pi_0)=\Pi_1$ where $\Pi_0=dx\wedge dy$ is the tangent Poisson lift of $(\nabla,h_0)$ and $\Pi_1=e^xdx\wedge dy$ is the tangent Poisson lift of $(\nabla,h_1)$.
\end{example}

\begin{example}[Different connections, same tensor]
	\label{ex:different-connections}
	On $\mathbb R$, consider
	\[
	\nabla^0_{\partial_x}\partial_x=0,\qquad
	\nabla^1_{\partial_x}\partial_x=\partial_x,\qquad
	h=\partial_x\otimes\partial_x.
	\]
	Both connections are flat and torsion-free. A diffeomorphism
	preserving $h$ must satisfy $(f')^2=1$, so $f(x)=\pm x+b$.
	For such a diffeomorphisms, 
	\[
	f_*(\nabla^0_{\partial_x}\partial_x)=0,
	\qquad
	\nabla^1_{f_*\partial_x}f_*\partial_x=\partial_x.
	\]
	The two KV manifolds are therefore not isomorphic. Their Poisson
	lifts nevertheless coincide. Indeed, the tangent Poisson lift of $(\nabla^0,h)$ is clearly 
	$\Pi^0=\partial_x\wedge\partial_y$ and we are going to prove that the tangent Poisson lift of $(\nabla^1,h)$ is $\Pi^1=\partial_x\wedge\partial_y$.
	Let's recall that if \(x\) is the standard coordinate on \(\mathbb R\), then a tangent vector
	\(u\in T_x\mathbb R\) has standard tangent coordinates \((x,y)\), where
	$y:=dx(u)$.
	We look for a coordinate \(t=t(x)\) adapted to the connection $\nabla^1$, i.e.\ such
	that the coordinate vector field \(\partial_t\) is parallel: $
	\nabla_{\partial_t}\partial_t=0$. We can see that the diffeomorphism 
	$$t:\mathbb R \to \mathbb R^*_+,\quad x\mapsto e^x$$
	is a global chart adapted to the connection $\nabla^1$ and the induced vector field $\partial_t$ is given by $\partial_t=e^{-x}\partial_x$ (which is clearly parallel vector field). The expression of the tensor $h$ in this chart is given by $h=t^2\partial_t\otimes\partial_t$. On the other hand, the tangent coordinates $(t,v)$ on $T\mathbb R\simeq\mathbb R^2$ induced from $t$ are given by:
	$$
	(t,v):u=y\ \frac{d}{dx}_{|_x}\in T\mathbb R\mapsto (t(x),dt_x(u))=\bigl(e^x,\;e^x y\bigr)\in \mathbb R^2
	$$
	hence $(t,v)=\bigl(e^x,\;e^x y\bigr)$. This implies that the tangent lift of $h$ is given by 
	$$\Pi^1=t^2\partial_t\wedge \partial_v=e^{2x}\left( (e^{-x}\partial_x)\wedge (e^{-x}\partial_y)\right)=\partial_x\wedge \partial_y=\Pi^0.
	$$
\end{example}

This leads to the following
question.

\medskip
\noindent\textbf{Question.}
\emph{Can the tangent Poisson structure be completed to a natural geometric structure that determines the underlying KV manifold up to isomorphism?}
\medskip

Our approach to addressing this question is inspired by Born geometry, which combines compatible complex and product structures with metric and symplectic data
\cite{FreidelLeighMinic2015,FRS2019}. Its pseudo-Riemannian formulation and
examples on Lie algebras have been studied in \cite{HKP2025,GP2025}.

In the pseudo-Hessian setting, Sakamoto \cite{Sakamoto2026} constructs a
strongly integrable Born structure on the tangent bundle and proves that it
determines the underlying pseudo-Hessian manifold up to isomorphism. Sakamoto
also describes its automorphism group in terms of pseudo-Hessian automorphisms
and parallel vector fields. This motivates a more general study adapted to the
context of KV manifolds, which is the primary objective of this article.

From this perspective, we introduce \emph{Born--Poisson} (BP) structures,
which combine a Poisson bivector \(\pi\) on a manifold \(M\) with a complex
product structure \((J,K)\) compatible with \(\pi\) (see
Section~\ref{sec:Complex product}. A BP isomorphism is a diffeomorphism
that preserves the Poisson bivectors and intertwines the complex
and product structures. In Section~\ref{subsec:canonical-lift},
we show that every KV manifold induces a canonical BP structure
$\mathcal B=(TM,\Pi,J,K)$ on its tangent bundle.
Our main results are as follows.

\begin{theorem}\label{thm:KV-reconstruction}
	Let $(M,\nabla,h)$ and $(M',\nabla',h')$ be KV manifolds.
	The following statements are equivalent:
	\begin{enumerate}
		\item $(M,\nabla,h)$ and $(M',\nabla',h')$ are isomorphic.
		\item Their canonical tangent BP manifolds $\mathcal B$ and
		$\mathcal B'$ are isomorphic.
	\end{enumerate}
\end{theorem}

\begin{theorem}\label{thm:KV-isomorphisms}
	Let $(M,\nabla,h)$ and $(M',\nabla',h')$ be KV manifolds.
	A diffeomorphism $F:TM\to TM'$ is a BP isomorphism if and only if
	\[
	F=\tau_Z\circ Tf,
	\]
	where $f:(M,\nabla,h)\to(M',\nabla',h')$ is a KV isomorphism,
	$Z$ is a $\nabla'$-parallel vector field, and
	$\tau_Z(v_{x'})=v_{x'}+Z_{x'}$. The pair $(f,Z)$ is unique.
	Consequently,
	\[
	\Aut_{\mathrm{BP}}(\mathcal B)
	\cong\mathfrak X_\nabla(M)\rtimes
	\Aut_{\mathrm{KV}}(M,\nabla,h),
	\]
	where $\mathfrak X_\nabla(M)$ is the additive group of parallel
	vector fields and the action is $f\cdot Z=f_*Z$.
\end{theorem}

When $h$ is non degenerate, $h^{-1}$ is a pseudo-Hessian metric and
this theorems recovers \cite[Theorem~1.2,Theorem~1.3]{Sakamoto2026}.

Let $(M,\nabla)$ be an affine manifold, and let $(J,K)$ denote
the canonical complex product structure on $TM$. Given a Poisson
bivector field $\pi$ on $TM$, the preceding construction naturally
raises the following question.

\medskip
\noindent\textbf{Question.}
\emph{Is the assumption that $(\pi,J,K)$ is a Born--Poisson structure on $TM$
	sufficient to ensure the existence of a KV bivector field $h$ on
	$(M,\nabla)$ whose tangent Poisson lift $\Pi^h$ coincides with $\pi$?
	If not, what additional condition characterizes such lifts?}
\medskip

The Born--Poisson condition alone is not sufficient, as the following
example shows.

\begin{example}
	Let $M=\mathbb{R}$ with its standard flat connection, and let $(x,y)$
	be the induced coordinates on $TM\simeq\mathbb{R}^2$. Consider the Poisson bivector field on $TM$ given by $\pi=(1+y^2)\partial_x\wedge\partial y$. For the
	canonical complex product structure $(J,K)$ on $TM$, one can easily see that $(\pi,J,K)$ is a
	Born--Poisson structure on $TM$. However, every KV bivector field on $M$ has the form $h=a(x) \partial x\otimes \partial x$
	and its tangent Poisson lift on $TM$ is
	\[
	\Pi=a(x)\  \partial x\wedge \partial y.
	\]
\end{example}

In response to this query, lets denote by $E$ the Euler vector field on $TM$, whose
flow is given by $\Phi_t(u)=e^t u,\, \forall u\in TM$. We have the following characterization.
\begin{theorem}\label{theo:homogeneous-Poisson-KV}
	Let $(M,\nabla)$ be an affine manifold, $J$ be its canonical
	complex structure on $TM$ and $\pi$ a bivector field on $TM$. The following conditions are equivalent:
	\begin{enumerate}
		\item[(i)] There exists a KV bivector $h$ on $(M,\nabla)$ whose tangent Poisson lift $\Pi^h$ coincides with
		$\pi$.
		\item[(ii)] The bivector $\pi$ is Poisson and satisfies
		\[
		\mathcal L_E\pi=-\pi,\qquad
		\pi(J^*\alpha,J^*\beta)=\pi(\alpha,\beta),
		\quad \forall\alpha,\beta\in\Omega^1(TM).
		\]
	\end{enumerate}
	Whenever any of these conditions hold, $h$ is unique.
\end{theorem}

The paper is organized as follows.
Section~\ref{sec:BP-structures} recalls complex product structures,
defines BP and symplectic-Born structures, and describes symplectic
leaves and BP structures on Lie algebra duals.
Section~\ref{sec:tangent-BP} constructs the canonical tangent BP
structure and we proves the homogeneity characterization.
Section~\ref{sec:reconstruction} proves the reconstruction and
automorphism theorems and identifies the Euler-preserving
isomorphisms. We applies these results to duals of
commutative associative algebras and computes both their BP and
Euler-preserving automorphism groups.

\medskip
\noindent\textbf{Notation.} For a contravariant
$2$-tensor $Q$, our convention is
$\beta(Q_\sharp\alpha)=Q(\alpha,\beta)$.

\section{Born--Poisson structures}\label{sec:BP-structures}
We  recall the complex and product structures underlying
Born--Poisson geometry (\cite{AS2005,Andrada2005,CFG1996,YanoAko1973}). Next, we introduce the notion of Born-Poisson structure. 
\begin{definition}[Poisson structure]
	A \emph{Poisson structure} on a smooth manifold \(M\) is a bivector field
	\(\pi\in\Gamma(\Lambda^2 TM)\) satisfying \([\pi,\pi]_{SN}=0\), where
	\([\cdot,\cdot]_{SN}\) is the Schouten--Nijenhuis bracket.
\end{definition}

We use the convention \(\beta(\pi_\sharp\alpha)=\pi(\alpha,\beta)\).
The characteristic distribution \(\operatorname{Im}\pi_\sharp\) is integrable
in the sense of singular foliations. Each maximal connected integral
manifold \(S\) carries a canonical symplectic form \(\omega_S\) given by
\[
\omega_S(\pi_\sharp\alpha,\pi_\sharp\beta)=\pi(\alpha,\beta).
\]
The submanifolds \(S\) are the symplectic leaves of \(M\)
(see \cite{Camille}).

\subsection{Complex product and Born--Poisson structures}\label{sec:Complex product}

\begin{itemize}
	\item A \emph{complex structure} on $M$ is an endomorphism
	$J\in\Gamma(\operatorname{End}(TM))$ satisfying
	$J^2=-\operatorname{id}_{TM}$ and $N_J=0$, where
	\[
	N_J(X,Y)
	=[JX,JY]-J[JX,Y]-J[X,JY]-[X,Y],
	\qquad X,Y\in\Gamma(TM).
	\]
	
	\item A \emph{product structure} on $M$ is an endomorphism
	$K\in\Gamma(\operatorname{End}(TM))$ satisfying
	$K^2=\operatorname{id}_{TM}$, $K\ne\pm\operatorname{id}_{TM}$,
	and $N_K=0$, where
	\[
	N_K(X,Y)
	=[KX,KY]-K[KX,Y]-K[X,KY]+[X,Y],
	\qquad X,Y\in\Gamma(TM).
	\]
	Its eigenbundles are
	\[
	T^+M=\ker(K-\operatorname{id}_{TM}),
	\qquad
	T^-M=\ker(K+\operatorname{id}_{TM}).
	\]
	The condition $N_K=0$ is equivalent to the involutivity of both
	eigenbundles. If they have the same rank, $K$ is called a
	\emph{para-complex structure}.
	
	\item A \emph{complex product structure} on $M$ is a pair
	$(J,K)$ consisting of a complex structure $J$ and a product
	structure $K$ such that
	\[
	KJ=-JK.
	\]
	The anticommutation relation implies that $J$ interchanges the
	eigenbundles of $K$. They therefore have the same rank, so $K$ is automatically
	para-complex. 
\end{itemize}

\begin{definition}[Born--Poisson structure]\label{def:Born-Poisson}
	A \emph{Born--Poisson} (BP) structure on $M$ is a triple
	$(\pi,J,K)$, where $\pi$ is a Poisson
	bivector and $(J,K)$ is a complex product structure satisfying
	\begin{equation}\label{eq:Pi-compatibility}
	\pi(J^*\alpha,J^*\beta)=\pi(\alpha,\beta),
	\qquad
	\pi(K^*\alpha,K^*\beta)=-\pi(\alpha,\beta),
	\end{equation}
	for all $\alpha,\beta\in\Omega^1(M)$.
	The quadruple $(M,\pi,J,K)$ is called a \emph{BP manifold}.
\end{definition}

The compatibility identities~\eqref{eq:Pi-compatibility} are
equivalently written as
\begin{equation}\label{eq:BP-anchor-compatibility}
J\circ\pi_\sharp=-\pi_\sharp\circ J^*,
\qquad
K\circ\pi_\sharp=\pi_\sharp\circ K^*.
\end{equation}
In particular, the characteristic distribution
$\operatorname{Im}\pi_\sharp$ is invariant under both $J$ and $K$.

\subsection{Symplectic-Born structures and symplectic leaves}

\begin{definition}[Symplectic-Born structure]\label{def:symplectic-Born}
	A \emph{symplectic-Born structure} on $M$ is a triple
	$(\omega,J,K)$, where $\omega$ is a symplectic $2$-form and
	$(J,K)$ is a complex product structure satisfying
	\begin{equation}\label{eq:pseudo-Born-compatibility}
	\omega(JX,JY)=\omega(X,Y),
	\qquad
	\omega(KX,KY)=-\omega(X,Y),
	\end{equation}
	for all $X,Y\in\Gamma(TM)$.
	The quadruple $(M,\omega,J,K)$ is called a
	\emph{symplectic-Born manifold}.
\end{definition}

\begin{remark}
	Definition~\ref{def:symplectic-Born} is precisely the symplectic formulation of an integrable
	Born manifold in the pseudo-Riemannian sense of \cite{HKP2025}.
	Indeed, the associated metrics are recovered from
	$(\omega,J,K)$ by
	\begin{equation}\label{eq:pseudo-Born-metrics}
	g(X,Y)=\omega(JX,Y),
	\qquad
	\eta(X,Y)=\omega(KX,Y),
	\end{equation}
	and the third endomorphism is $JK$. Thus the symplectic and
	metric formulations describe the same geometric data.
\end{remark}

If $(M,\pi,J,K)$ is a BP manifold and $\pi$ is nondegenerate,
its associated symplectic form $\omega=\pi^{-1}$, characterized by
\[
\omega(\pi_\sharp\alpha,\pi_\sharp\beta)=\pi(\alpha,\beta),
\]
makes
$(M,\omega,J,K)$ a symplectic-Born manifold. The same construction applies on each positive-dimensional leaf.

\begin{proposition}\label{prop:symplectic-leaves}
	Let $(M,\pi,J,K)$ be a BP manifold, and let $S$ be a symplectic leaf of $(M,\pi)$.
	If $\dim S>0$, then $(S,\omega_S,J_S,K_S)$ is a symplectic-Born
	manifold, where $\omega_S$ is the symplectic form of $S$ and
	\[
	J_S=J|_{TS},
	\qquad
	K_S=K|_{TS}.
	\]
\end{proposition}

\begin{proof}
	Since $TS=\operatorname{Im}\pi_\sharp|_S$, the identities \eqref{eq:BP-anchor-compatibility}
	show that $TS$ is invariant under $J$ and $K$.
	Their restrictions satisfy
	\[
	J_S^2=-\operatorname{id}_{TS},
	\qquad
	K_S^2=\operatorname{id}_{TS},
	\qquad
	K_SJ_S=-J_SK_S.
	\]
	Moreover, the Nijenhuis tensors of $J_S$ and $K_S$ are the
	restrictions of the ambient Nijenhuis tensors and therefore
	vanish. Hence $(J_S,K_S)$ is a complex product structure.
	
	Write $X=\pi_\sharp\alpha$ and
	$Y=\pi_\sharp\beta$. Using the identities \eqref{eq:Pi-compatibility}, we obtain
	\[
	\begin{aligned}
	\omega_S(J_SX,J_SY)
	&=\pi(J^*\alpha,J^*\beta)
	=\pi(\alpha,\beta)
	=\omega_S(X,Y),\\
	\omega_S(K_SX,K_SY)
	&=\pi(K^*\alpha,K^*\beta)
	=-\pi(\alpha,\beta)
	=-\omega_S(X,Y).
	\end{aligned}
	\]
	Thus $(S,\omega_S,J_S,K_S)$ is a symplectic-Born manifold.
\end{proof}

\section{The canonical tangent Born--Poisson structure}
\label{sec:tangent-BP}

\subsection{Horizontal--vertical splitting and the Poisson lift}
\label{subsec:canonical-lift}

Let $(M,\nabla)$ be an affine manifold, and let $p:TM\to M$ be the canonical projection. The connection $\nabla$ determines a splitting
\[
T(TM)=H\oplus V,\qquad V=\ker(Tp).
\]
For $X\in\mathfrak X(M)$, its horizontal and vertical lifts are
defined by
\[
X^h_u=(T_up|_{H_u})^{-1}(X_{p(u)}),\qquad
X^v_u=\left.\frac{d}{dt}\right|_{t=0}
\bigl(u+tX_{p(u)}\bigr).
\]
Using these lifts, define $J,K\in\Gamma(\End(T(TM)))$ by
\begin{equation}\label{eq:canonical-AB}
	\begin{cases}
		J(X^h)=X^v,\\ J(X^v)=-X^h,
	\end{cases}
	\qquad
	\begin{cases}
		K(X^h)=X^h,\\ K(X^v)=-X^v.
	\end{cases}
\end{equation}
Then $J^2=-\id$, $K^2=\id$, and $KJ=-JK$.
In induced affine coordinates, both endomorphisms have constant
coefficients. Their Nijenhuis tensors therefore vanish, so $(J,K)$
is a complex product structure. 

For $\alpha\in\Omega^1(M)$, its horizontal and vertical lifts are
the $1$-forms characterized by
\[
\begin{cases}
	\alpha^h(X^h)=\alpha(X)\circ p,\\ \alpha^h(X^v)=0,
\end{cases}
\qquad
\begin{cases}
	\alpha^v(X^h)=0,\\ \alpha^v(X^v)=\alpha(X)\circ p.
\end{cases}
\]
A symmetric bivector $h$ on $M$ determines a skew-symmetric bivector field $\Pi$ on
$TM$ by
\begin{equation}\label{eq:lift-pairings}
	\Pi(\alpha^h,\beta^h)=\Pi(\alpha^v,\beta^v)=0,
	\qquad
	\Pi(\alpha^h,\beta^v)=h(\alpha,\beta)\circ p.
\end{equation}
In induced affine coordinates,
\begin{equation}\label{eq:canonical-Pi}
	\Pi=\sum_{i,j=1}^n h_{ij}
	\partial_{x_i}\wedge\partial_{y_j}.
\end{equation}

\begin{theorem}\label{thm:construction}
	Let $(M,\nabla,h)$ be a KV manifold. Then $(TM,\Pi,J,K)$ is
	a BP manifold.
\end{theorem}

\begin{proof}
	By \cite[Theorem~3.3]{ABB2020}, $\Pi$ is Poisson. The complex
	product structure was established above. The dual endomorphisms
	satisfy
	\[
	J^*\alpha^h=-\alpha^v,\qquad J^*\alpha^v=\alpha^h,
	\qquad
	K^*\alpha^h=\alpha^h,\qquad K^*\alpha^v=-\alpha^v.
	\]
	Thus skew-symmetry of $\Pi$ and symmetry of $h$ give
	\[
	\begin{aligned}
		\Pi(J^*\alpha^h,J^*\beta^v)
		&=\Pi(-\alpha^v,\beta^h)
		=h(\beta,\alpha)\circ p
		=\Pi(\alpha^h,\beta^v),\\
		\Pi(K^*\alpha^h,K^*\beta^v)
		&=\Pi(\alpha^h,-\beta^v)
		=-\Pi(\alpha^h,\beta^v).
	\end{aligned}
	\]
	For pairs of horizontal forms or pairs of vertical forms, the
	identities follow from the vanishing pairings in
	\eqref{eq:lift-pairings}. So both compatibility identities \ref{eq:pseudo-Born-compatibility} hold for all $1$-forms.
\end{proof}

\subsection{Homogeneity and characterization of KV lifts}
\label{subsec:homogeneity}

To prove the characterization stated in
Theorem~\ref{theo:homogeneous-Poisson-KV} we need the following lemma.
\begin{lemma} Let $(M,\nabla)$ be an affine manifold. For any function $f\in C^\infty(M)$, any vector field $X\in\mathfrak{X}(M)$, and any one form $\alpha\in\Omega^1(M)$, 
	\begin{equation}\label{eq:Euler-lifts} E(f\circ p)=0,\quad 
		[E,X^h]=0,\quad [E,X^v]=-X^v,
		\quad
		\mathcal L_E\alpha^h=0,\quad
		\mathcal L_E\alpha^v=\alpha^v
	\end{equation}	
\end{lemma}
\begin{proof}
	Let $\Phi_t(u)=e^t u$ be the flow of $E$. Since
	$p\circ\Phi_t=p$, we have
	\[
	(f\circ p)\circ\Phi_t=f\circ p.
	\]
	Differentiating at $t=0$ gives $E(f\circ p)=0$.
	
	Since parallel transport is linear, fibre dilations preserve
	horizontal lifts. Moreover, by the definition of vertical lifts,
	\[
	T_u\Phi_t(X^h_u)=X^h_{\Phi_t(u)},
	\qquad
	T_u\Phi_t(X^v_u)=e^tX^v_{\Phi_t(u)}.
	\]
	Thus the pullbacks of the lifted vector fields satisfy
	\[
	\Phi_t^*X^h=X^h,
	\qquad
	\Phi_t^*X^v=e^{-t}X^v.
	\]
	Differentiating at $t=0$ yields
	\[
	[E,X^h]=0,
	\qquad
	[E,X^v]=-X^v.
	\]
	Similarly, we obtain
	\[
	\mathcal L_E\alpha^h=0,
	\qquad
	\mathcal L_E\alpha^v=\alpha^v.
	\]
\end{proof}

\begin{proof}[Proof of Theorem~\ref{theo:homogeneous-Poisson-KV}]
	Suppose first that $\pi=\Pi$ is the tangent Poisson lift of a
	KV bivector $h$. By Theorem~\ref{thm:construction}, $\pi$ is Poisson
	and compatible with $J$. Equation~\eqref{eq:Euler-lifts} gives
	$\mathcal L_E\pi=-\pi$. Thus $(i)$ implies $(ii)$.
	
	Conversely, suppose that $(ii)$ holds. Recall that
	\[
	J^*\alpha^h=-\alpha^v,
	\qquad
	J^*\alpha^v=\alpha^h.
	\]
	Compatibility with $J$ therefore gives
	\begin{equation}\label{eq:J-lift-pairings}
		\pi(\alpha^v,\beta^v)=\pi(\alpha^h,\beta^h),
		\qquad
		\pi(\alpha^h,\beta^v)=\pi(\beta^h,\alpha^v)
	\end{equation}
	for all $\alpha,\beta\in\Omega^1(M)$.
	
	Using $\mathcal L_E\pi=-\pi$ and \eqref{eq:Euler-lifts}, we obtain
	\[
	E\bigl(\pi(\alpha^h,\beta^h)\bigr)
	=-\pi(\alpha^h,\beta^h),
	\qquad
	E\bigl(\pi(\alpha^v,\beta^v)\bigr)
	=\pi(\alpha^v,\beta^v).
	\]
	From the first identity in \eqref{eq:J-lift-pairings}, we get
	\begin{equation*}\label{eq:pure-lift-vanishing}
		\pi(\alpha^h,\beta^h)=\pi(\alpha^v,\beta^v)=0.
	\end{equation*}
	The same calculation yields
	\[
	E\bigl(\pi(\alpha^h,\beta^v)\bigr)=0.
	\]
	Thus, for any $u\in TM$,
	\[
	\pi_{tu}(\alpha^h_{tu},\beta^v_{tu})
	=\pi_u(\alpha^h_u,\beta^v_u),
	\qquad \forall t>0,
	\]
	letting $t\to0^+$ gives
	\[
	\pi_u(\alpha^h_u,\beta^v_u)
	=\pi_{0_x}(\alpha^h_{0_x},\beta^v_{0_x}),
	\qquad x=p(u).
	\]
	
	Define a tensor $h$ on $M$ by
	\begin{equation}\label{eq:recover-h-zero}
		h_x(\alpha_x,\beta_x)
		=\pi_{0_x}(\alpha^h_{0_x},\beta^v_{0_x}).
	\end{equation}
	The right-hand side depends only on $\alpha_x$ and $\beta_x$.
	It defines a smooth bilinear tensor, which is symmetric by
	\eqref{eq:J-lift-pairings}. Moreover,
	\[
	\pi(\alpha^h,\beta^h)=\pi(\alpha^v,\beta^v)=0,
	\qquad
	\pi(\alpha^h,\beta^v)=h(\alpha,\beta)\circ p.
	\]
 Hence $\pi=\Pi^h$.
	
	From \cite[Theorem~3.3]{ABB2020} It follows that $h$ is a KV bivector field. Its uniqueness follows immediately from
	\eqref{eq:recover-h-zero}.
\end{proof}

\section{Reconstruction and automorphisms}
\label{sec:reconstruction}

\subsection{Rigidity of the canonical complex product structure}
\label{subsec:rigidity}

Let $(M,\nabla)$ and $(M',\nabla')$ be affine manifolds, with canonical 
projections $p:TM\rightarrow M$,$p':TM'\rightarrow M'$ and complex product structure $(J,K)$ and $(J',K')$ given by \ref{eq:canonical-AB}. For $Z\in\mathfrak X(M')$, write
\[
\tau_Z:TM'\to TM',\qquad \tau_Z(v_{x'})=v_{x'}+Z_{x'}.
\]
The proof of the following lemma is adapted from  Proposition~4.23 in \cite{Sakamoto2026} to the case of flat manifolds.

\begin{lemma}\label{lem:lift-rigidity}
	Let $F:TM\to TM'$ be a diffeomorphism. Then
	
	\[ TF\circ J=J'\circ TF, \text{ and }TF\circ K=K'\circ TF\]
	if and only if, there exists an affine diffeomorphism $f:(M,\nabla)\to(M',\nabla')$  
	and a $\nabla'$-parallel vector field $Z\in\mathfrak{X}(M')$ such that
	\[
	F=\tau_Z\circ Tf.
	\]
	In particular, the pair $(f,Z)$ is unique.
\end{lemma}

\begin{proof}
	
	\medskip
	\noindent
	\textbf{Step 1. The definition of $f$ and $Z$.}
	Assume that
	\[
	TF\circ J=J'\circ TF,
	\qquad
	TF\circ K=K'\circ TF.
	\]	
	Since $H$ and $V$ are the $+1$ and $-1$ eigendistributions
	of $K$, respectively, the relation $TF\circ K=K'\circ TF$,
	implies
	\[
	TF(H)=H',
	\qquad
	TF(V)=V'.
	\]	
	In particular, for every $u\in TM$ we have $T_u(p'\circ F)(V_u)=\{0\}$. Thus the map $p'\circ F$
	is constant on $T_xM$. Consequently, there exists a unique smooth map
	$f:M\to M'$ such that $	p'\circ F=f\circ p$. Explicitly, $f(x)=p'\bigl(F(0_x)\bigr)$. Thus, $f$ is a diffeomorphism.
	
Fix $x\in M$. Let $u,v\in T_xM$ and denote by $F_x=F|_{T_xM}:T_xM\rightarrow T_{f(x)}M'$.
		We claim that 
		\begin{equation}\label{eq:fibre-derivative}
			T_v(F_x)=T_xf.
		\end{equation}
	Indeed, using the relation $TF\circ J=J'\circ TF$, we obtain
	\begin{equation*}\label{eq:horizontal-identity}
		T_vF(u_v^{\mathrm h})
		=
		\bigl(T_v(F_x)(u)\bigr)_{F(v)}^{\mathrm h}.
	\end{equation*}
	Applying $T_{F(v)}p'$ to both sides gives
	\[
	T_{F(v)}p'\bigl(T_vF(u_v^{\mathrm h})\bigr)
	=
	T_v(F_x)(u).
	\]
		On the other hand, differentiating
	$p'\circ F=f\circ p$, we have $T_{F(v)}p'\circ T_vF
	=
	T_xf\circ T_vp$. Hence
	\[
	T_v(F_x)(u)=T_xf(u).
	\]
	For any $v\in T_xM$, consider the curve $	c(t)=F_x(tv),
	\,\, 0\leq t\leq 1.$
	By \eqref{eq:fibre-derivative},
	\[
	c'(t)=T_{tv}(F_x)(v)=T_xf(v).
	\]
	Integrating from $0$ to $1$, we obtain $	F_x(v)-F_x(0_x)=T_xf(v)$.
	Hence
	\begin{equation}\label{eq:fibre-affine}
		F(v)=T_xf(v)+F(0_x).
	\end{equation}
		Define a vector field $Z$ on $M'$ by
	\begin{equation}\label{eq:Z-definition}
		Z_{x'}
		=
		F\bigl(0_{f^{-1}(x')}\bigr),
		\qquad x'\in M'.
	\end{equation}
 By construction, $F(0_x)=Z_{f(x)}$.
	Thus \eqref{eq:fibre-affine} becomes
	\begin{equation}\label{eq:F-decomposition}
		F(v)=T_xf(v)+Z_{f(x)},
		\qquad v\in T_xM.
	\end{equation}
	This means that $F=\tau_Z\circ Tf$.
	
\medskip
\noindent
\textbf{Step 2. $Z$ is $\nabla'$-parallel.}	
Let $\gamma:I\to M$ be a smooth curve and
consider the zero vector field along $\gamma$, $U(t)=0_{\gamma(t)}$.
Since $\nabla_{\dot\gamma}U=0$, the curve $U$ is
horizontal in $TM$. As $TF(H)=H'$, its image
$F\circ U$ is horizontal in $TM'$. On the other hand, by \eqref{eq:F-decomposition},
\[
F(U(t))=F(0_{\gamma(t)})=Z_{f(\gamma(t))}.
\]
Consequently, $Z$ is parallel along $f\circ\gamma$,
that is,
\begin{equation}\label{eq:Z-parallel}
	\nabla'_{T_{\gamma(t)}f(\dot\gamma(t))}Z=0.
\end{equation}
Now fix $x'\in M'$ and $u'\in T_{x'}M'$, and set
$x=f^{-1}(x')$ and $u=(T_xf)^{-1}(u')$.
Choose a smooth curve $\gamma:I\rightarrow M$ satisfying $\gamma(0)=x,\,\, \dot\gamma(0)=u$.
Evaluating the identity \eqref{eq:Z-parallel} at $t=0$
gives $(\nabla'_{u'}Z)_{x'}=0$. Thus, $\nabla'Z=0$.

	\medskip
	\noindent
	\textbf{Step 3. $f$ is affine.} Let $\gamma:I\to M$ be a smooth curve and  $U$ a  $\nabla$-parallel vector field
	along $\gamma$. We are going to prove that $Tf\circ U$ is $\nabla'$-parallel vector field
	along $\alpha:=f\circ\gamma$.
	Since $U$ is horizontal and $TF(H)=H'$, the curve $F\circ U$ is horizontal in $TM'$.
	Using \eqref{eq:F-decomposition}, we can write
	\[
	F(U(t))
	=
	T_{\gamma(t)}f(U(t))
	+
	Z_{\alpha(t)}.
	\]
	Taking the covariant derivative along $\alpha$,
	we obtain
	\[ 0=\nabla'_{\dot\alpha}(F\circ U)=\nabla'_{\dot\alpha}
	\bigl(Tf\circ U\bigr)
	+\nabla'_{\dot\alpha}Z.
	\]
	Since $\nabla'Z=0$, this reduces to $\nabla'_{\dot\alpha}
	\bigl(Tf\circ U\bigr)=0$.
	
	\medskip
	\noindent
	\textbf{Step 5. The converse.}
	Conversely, suppose that $f:(M,\nabla)\longrightarrow(M',\nabla')$
	is an affine diffeomorphism and that
	$Z\in\mathfrak X(M')$ satisfies $\nabla'Z=0$. 	
	Define
	\[
	F(v)=T_xf(v)+Z_{f(x)},
	\qquad v\in T_xM.
	\]
	Obviously $F$ is a diffeomorphism. We show that $TF$ preserves both the horizontal
	and vertical lifts. Let $x\in M$ and $u,v\in T_xM$
	\[
	\begin{aligned}
		T_vF(u_v^{\mathrm v})
		&=
		\left.\frac{d}{dt}\right|_{t=0}F(v+tu)\\
		&=
		\left.\frac{d}{dt}\right|_{t=0}
		\left(
		T_xf(v+tu)+Z_{f(x)}
		\right)\\
		&=
		(T_xf(u))_{F(v)}^{\mathrm v}.
	\end{aligned}
	\]
	Therefore
	\begin{equation}\label{eq:converse-vertical}
		T_vF(u_v^{\mathrm v})
		=
		(T_xf(u))_{F(v)}^{\mathrm v}.
	\end{equation}
	Next, choose a smooth curve $\gamma:I\rightarrow M$ satisfying, $\gamma(0)=x,
	\,\,
	\dot\gamma(0)=u,$
	and let $U$ be the $\nabla$-parallel vector field
	along $\gamma$ with $U(0)=v$. By the definition of the horizontal lift,
	\[
	u_v^{\mathrm h}
	=
	\left.\frac{d}{dt}\right|_{t=0}U(t).
	\]
		Since $f$ is affine, $Tf\circ U$ is
	$\nabla'$-parallel along $f\circ\gamma$.
	Since $\nabla'Z=0$, the vector field
	\[
	F(U(t))
	=
	Tf(U(t))+Z_{f(\gamma(t))}
	\]
	is also $\nabla'$-parallel along $f\circ\gamma$. Consequently,
	\[ 	(T_xf(u))_{F(v)}^{\mathrm h}=\left.\frac{d}{dt}\right|_{t=0}F(U(t))=T_vF(u_v^{\mathrm h}).
	\]
	Thus
	\begin{equation}\label{eq:converse-horizontal}
		T_vF(u_v^{\mathrm h})
		=
		(T_xf(u))_{F(v)}^{\mathrm h}.
	\end{equation}
	Using \eqref{eq:converse-vertical} and
	\eqref{eq:converse-horizontal}, we obtain
	\[
	\begin{aligned}
		T_vF\bigl(J(u_v^{\mathrm h})\bigr)
		&=T_vF(u_v^{\mathrm v})\\
		&=(T_xf(u))_{F(v)}^{\mathrm v}\\
		&=J'\bigl((T_xf(u))_{F(v)}^{\mathrm h}\bigr)\\
		&=J'\bigl(T_vF(u_v^{\mathrm h})\bigr).
	\end{aligned}
	\]
	Similarly,
	\[
	\begin{aligned}
		T_vF\bigl(J(u_v^{\mathrm v})\bigr)
		&=-T_vF(u_v^{\mathrm h})\\
		&=-(T_xf(u))_{F(v)}^{\mathrm h}\\
		&=J'\bigl((T_xf(u))_{F(v)}^{\mathrm v}\bigr)\\
		&=J'\bigl(T_vF(u_v^{\mathrm v})\bigr).
	\end{aligned}
	\]
Hence, $TF\circ J=J'\circ TF$. In the same way, we get that $TF\circ K=K'\circ TF$.
\end{proof}

\subsection{Poisson naturality and the main theorems}
\label{subsec:main-proofs}

\begin{lemma}\label{prop:Poisson-naturality}
	Let $f:(M,\nabla,h)\to(M',\nabla',h')$ be an affine diffeomorphism between KV manifolds,  and let $Z$ be $\nabla'$-parallel  vector field on $M'$. Denote by $\Pi$ and $\Pi'$ the Poisson tensor on $TM$ and $TM'$, associated to $(\nabla,h)$ and $(\nabla',h')$, respectively. 
	For $F=\tau_Z\circ Tf$, one has
	\[
	F_*\Pi=\Pi'\quad\Longleftrightarrow\quad f_*h=h'.
	\]
\end{lemma}

\begin{proof}
From equations \eqref{eq:converse-vertical} and
	\eqref{eq:converse-horizontal} one can see that, for every
	$\alpha'\in T_{f(x)}^*M'$,
	\begin{equation}\label{eq:F-pullback-covectors}
		\begin{aligned}
			(T_vF)^*(\alpha_{v'}^{\prime\mathrm h})
			&=(f^*\alpha')_v^{\mathrm h},\\
			(T_vF)^*(\alpha_{v'}^{\prime\mathrm v})
			&=(f^*\alpha')_v^{\mathrm v},
		\end{aligned}
	\end{equation}
	where $f^*\alpha'=\alpha'\circ T_xf$.
	
	Let $\alpha',\beta'\in T_{f(x)}^*M'$. 
	From the definition of $\Pi$ and \eqref{eq:F-pullback-covectors}
	\[
	\begin{aligned}
		(F_*\Pi)_{v'}
		(\alpha_{v'}^{\prime\mathrm v},
		\beta_{v'}^{\prime\mathrm h})
		&=
		\Pi_v\bigl(
		(T_vF)^*\alpha_{v'}^{\prime\mathrm v},
		(T_vF)^*\beta_{v'}^{\prime\mathrm h}
		\bigr)\\
		&=\Pi_v\bigl(
		(f^*\alpha')_v^{\mathrm v},
		(f^*\beta')_v^{\mathrm h}
		\bigr)\\
		&=h_x(f^*\alpha',f^*\beta')\\
		&=(f_*h)_{f(x)}(\alpha',\beta').
	\end{aligned}
	\]
	Consequently,
	\begin{equation*}\label{eq:Pi-naturality-mixed}
			(F_*\Pi)_{v'}
			(\alpha_{v'}^{\prime\mathrm v},
			\beta_{v'}^{\prime\mathrm h})
			=
			(f_*h)_{f(x)}(\alpha',\beta').
	\end{equation*}
		Similarly, we obtain
	\[
	\begin{aligned}
		(F_*\Pi)_{v'}
		(\alpha_{v'}^{\prime\mathrm h},
		\beta_{v'}^{\prime\mathrm h})&=0=	\Pi'_{v'},
		(\alpha_{v'}^{\prime\mathrm h},
		\beta_{v'}^{\prime\mathrm h}),\\
		(F_*\Pi)_{v'}
		(\alpha_{v'}^{\prime\mathrm v},
		\beta_{v'}^{\prime\mathrm v})&=0=\Pi'_{v'}
		(\alpha_{v'}^{\prime\mathrm v},
		\beta_{v'}^{\prime\mathrm v}).
	\end{aligned}
	\]
And the above equivalence follows.
\end{proof}

\begin{proof}[Proof of Theorem~\ref{thm:KV-reconstruction}]
	If $f$ is a KV isomorphism, Lemma~\ref{lem:lift-rigidity} with
	$Z=0$ shows that $Tf$ preserves the canonical complex product $(J,K)$.
	Lemma~\ref{prop:Poisson-naturality} shows that it preserves the
	Poisson bivectors. Thus $Tf$ is a BP isomorphism.
	
	Conversely, a BP isomorphism $F$ preserves $(J,K)$, so
	Lemma~\ref{lem:lift-rigidity} gives $F=\tau_Z\circ Tf$, with
	$f$ affine and $Z$ parallel. Preservation of the Poisson bivectors
	then gives $f_*h=h'$ by Lemma~\ref{prop:Poisson-naturality}.
	Thus $f$ is a KV isomorphism.
\end{proof}

\begin{proof}[Proof of Theorem~\ref{thm:KV-isomorphisms}]
	Lemmas~\ref{lem:lift-rigidity} and~\ref{prop:Poisson-naturality}
	give the stated classification and uniqueness. Every KV
	automorphism is affine and sends parallel vector fields to parallel
	vector fields. It therefore acts on $\mathfrak X_\nabla(M)$ by
	pushforward. Hence we get a bijection
	\[
	\begin{aligned}
		\Phi:\mathfrak X_\nabla(M)\rtimes
		\Aut_{\mathrm{KV}}(M,\nabla,h)
		&\longrightarrow\Aut_{\mathrm{BP}}(\mathcal B),\\
		(Z,f)&\longmapsto\tau_Z\circ Tf.
	\end{aligned}
	\]
	Since $Tf\circ\tau_W=\tau_{f_*W}\circ Tf$,
	\[
	(\tau_Z\circ Tf)\circ(\tau_W\circ Tg)
	=\tau_{Z+f_*W}\circ T(f\circ g).
	\]
	This is the semidirect-product law
	\begin{equation}\label{eq:aut-group-law}
		(Z,f)(W,g)=(Z+f_*W,f\circ g),
	\end{equation}
	so $\Phi$ is a group isomorphism.
\end{proof}

\subsection{Euler-preserving isomorphisms}
\label{subsec:Euler-isomorphisms}
A diffeomorphism $F:TM\to TM$ preserve $(\Pi,E,J)$ means that $F_*\Pi=\Pi$, $F_*E=E$, and
$TF\circ J=J\circ TF$.

\begin{proposition}\label{cor:Euler-isomorphisms}
	Let $(M,\nabla,h)$ be a KV manifold. Let $F:TM\to TM$ be a
	a diffeomorphism. Then, $F$ preserve $(\Pi,E,J)$ if and only
	if, $F=Tf$, where  $f:(M,\nabla,h)\rightarrow (M,\nabla,h)$ is a KV isomorphism. In other words,
	\begin{equation*}\label{eq:Euler-aut-group}
		\Aut(TM,\Pi,E,J)\cong\Aut_{\mathrm{KV}}(M,\nabla,h).
	\end{equation*}
\end{proposition}

\begin{proof}	Assume first that $F$ preserves $(\Pi,E,J)$. Let $K:=J\circ\mathcal L_EJ$, obviously $K$ is the product structure on $TM$ associated to $\nabla$. Since $F_*E=E$ and $TF\circ J=J\circ TF$, for every
	$X\in\mathfrak X(TM)$, we have
	\[
	\begin{aligned}
		F_*\big((\mathcal L_EJ)(X)\big)
		&=F_*\big([E,JX]-J[E,X]\big)\\
		&=[F_*E,F_*(JX)]-J(F_*[E,X])\\
		&=[E,J(F_*X)]-J[E,F_*X]\\
		&=(\mathcal L_EJ)(F_*X).
	\end{aligned}
	\]
	Thus
	\[
	TF\circ\mathcal L_EJ
	=\mathcal L_EJ\circ TF.
	\]
	Hence, we obtain
	\[
	\begin{aligned}
		TF\circ K
		&=TF\circ J\circ\mathcal L_EJ\\
		&=J\circ TF\circ\mathcal L_EJ\\
		&=J\circ\mathcal L_EJ\circ TF\\
		&=K\circ TF.
	\end{aligned}
	\]
		Thus $F$ is a BP isomorphism. Theorem~\ref{thm:KV-isomorphisms} implies that
	there exist a KV isomorphism $	f:(M,\nabla)\longrightarrow(M,\nabla)$ and a $\nabla$-parallel vector field
	$Z\in\mathfrak X(M)$ such that for every $v\in T_xM$,
	\begin{equation*}\label{eq:Euler-F-decomposition}
		F(v)=T_xf(v)+Z_{f(x)}.
	\end{equation*}
	Since $F_*E=E$, we have
	\[
	E_{F(0_x)}
	=
	T_{0_x}F(E_{0_x})
	=0.
	\]	
	Thus $F(0_x)=0_{f(x)}$. Hence
	\[
	0_x=F(0_x)
	=
	T_xf(0_x)+Z_{f(x)}
	=
	Z_{f(x)}.
	\]	
	Consequently, $F=Tf$
	
	Conversely, let $f:(M,\nabla,h)\longrightarrow(M,\nabla,h)$
	be a KV isomorphism, and set $F=Tf$. Theorem~\ref{thm:KV-isomorphisms} implies that $F$ preserves $(\Pi,J,K)$.
	It remains to show that $F$ preserves $E$. For every $v\in T_xM$, we have
	\[
	\begin{aligned}
		T_vF(E_v)
		&=
		T_v(Tf)
		\left(
		\left.\frac{d}{dt}\right|_{t=0}e^tv
		\right)\\
		&=
		\left.\frac{d}{dt}\right|_{t=0}
		Tf(e^tv)\\
		&=
		\left.\frac{d}{dt}\right|_{t=0}
		e^tTf(v)\\
		&=
		E_{F(v)}.
	\end{aligned}
	\]
	Thus $F_*E=E$.
\end{proof}

\subsection{Example}

Let $(\mathcal A,\cdot)$ be a finite-dimensional  real
commutative associative algebra. Its dual $\mathcal A^*$ has the
standard flat connection $\nabla^0$ and the symmetric bivector
\begin{equation}\label{eq:algebra-KV}
	h_\xi(a,b)=\xi(a\cdot b),
	\qquad \xi\in\mathcal A^*,\quad
	a,b\in T_\xi^*\mathcal A^*\simeq\mathcal A.
\end{equation}
These data $(\nabla^0,h)$ define a KV structure on $\mathcal A^*$
\cite[Proposition~4.2]{ABB2020}.

Identify $T\mathcal A^*$ with $\mathcal A^*\oplus\mathcal A^*$,
write a point as $(\xi,\nu)$, and identify its cotangent space with
$\mathcal A\oplus\mathcal A$. The tangent Poisson bivector is
\begin{equation}\label{eq:algebra-Pi-intrinsic}
	\Pi_{(\xi,\nu)}\bigl((a,b),(c,d)\bigr)
	=\xi(a\cdot d-b\cdot c).
\end{equation} Equivalently, it
is the linear Lie--Poisson bivector on
$(\mathcal A\oplus\mathcal A)^*$ associated with Lie brackets
\[
[(a,b),(c,d)]=(a\cdot d-b\cdot c,0),
\]
 The complex product structure associated to $\nabla^0$ is given by 
\begin{equation}\label{eq:algebra-AJ}
	J(u,v)=(-v,u),\qquad K(u,v)=(u,-v).
\end{equation}
Theorem~\ref{thm:construction} gives that 
$\mathcal B=(T\mathcal A^*,\Pi,J,K)$ is a BP manifold.

Set
\[
\mathcal A^2=\operatorname{span}\{a\cdot b:a,b\in\mathcal A\},
\qquad
(\mathcal A^2)^\circ
=\{\alpha\in\mathcal A^*: \alpha(\mathcal A^2)=0\}.
\]
The group $\Aut_{\mathrm{alg}}(\mathcal A,\cdot)$ acts on $\mathcal A^*$ through
the contragredient representation
\[
\rho(\varphi)=(\varphi^{-1})^*.
\]
And the KV automorphism group is given by
\[ \Aut_{\mathrm{KV}}(\mathcal A^*,\nabla^0,h)=(\mathcal A^2)^\circ\rtimes_\rho
\Aut_{\mathrm{alg}}(\mathcal A,\cdot), \]
for the proof see \cite[Proposition 3.3]{ABB2023}. 

\begin{proposition}\label{prop:algebra-aut-BP}$\ $
	\begin{itemize}
		\item[$(i)$] Every $F\in Aut_{\mathrm{BP}}(\mathcal B)$ has the unique form
			\begin{equation}\label{eq:algebra-BP-map}
			F_{\alpha,\beta,\varphi}(\xi,\nu)
			=\bigl(\rho(\varphi)\xi+\alpha,
			\rho(\varphi)\nu+\beta\bigr),
		\end{equation}
		where $\alpha\in(\mathcal A^2)^\circ$, $\beta\in\mathcal A^*$,
		and $\varphi\in\Aut_{\mathrm{alg}}(\mathcal A,\cdot)$. Consequently,
		\begin{equation}\label{eq:algebra-aut-BP}
			\Aut_{\mathrm{BP}}(\mathcal B)
			\cong\bigl((\mathcal A^2)^\circ\oplus\mathcal A^*\bigr)
			\rtimes_\rho\Aut_{\mathrm{alg}}(\mathcal A,\cdot),
		\end{equation}
		\item[$(ii)$] The automorphism group that preserves $(\Pi,E,J)$ is given by \begin{equation}\label{eq:algebra-Euler-aut}
			\Aut(T\mathcal A^*,\Pi,E,J)
			\cong(\mathcal A^2)^\circ\rtimes_\rho
			\Aut_{\mathrm{alg}}(\mathcal A,\cdot).
		\end{equation}
	\item[$(iii)$]	If $\mathcal A$ is unital, then
		\begin{equation}\label{eq:unital-aut-BP}
			\begin{aligned}
				\Aut_{\mathrm{BP}}(\mathcal B)
				&\cong\mathcal A^*\rtimes_\rho\Aut_{\mathrm{alg}}(\mathcal A,\cdot),\\
				\Aut(T\mathcal A^*,\Pi,E,J)
				&\cong\Aut_{\mathrm{alg}}(\mathcal A,\cdot).
			\end{aligned}
		\end{equation}
	\end{itemize}
\end{proposition}

\begin{proof}$\ $
	\begin{itemize}
		\item[$(i)$] The parallel vector fields of $\nabla^0$ are the constant vector
		fields, naturally identified with $\mathcal A^*$. Their fibre
		translations are $(\xi,\nu)\mapsto(\xi,\nu+\beta)$. Hence
		Theorem~\ref{thm:KV-isomorphisms} gives
		\eqref{eq:algebra-BP-map}. Composition is
		\[
		F_{\alpha,\beta,\varphi}\circ F_{\alpha',\beta',\psi}
		=F_{\alpha+\rho(\varphi)\alpha',\,
			\beta+\rho(\varphi)\beta',\,\varphi\psi},
		\]
		which proves \eqref{eq:algebra-aut-BP}.
		\item[$(ii)$] Is a direct consequences of Proposition~\ref{cor:Euler-isomorphisms}.
	\item[$(iii)$] Follows from the fact that  $\mathcal A^2=\mathcal A$.
	\end{itemize}
\end{proof}

\textbf{Conflicts of interest}. The authors declare that they have no competing interests.

\textbf{Funding.} The authors received no financial support for this work.

\textbf{Data availability.} Not applicable to this work.


\bibliographystyle{amsplain}

\end{document}